\documentclass[11pt,a4paper,fleqn]{article}

\usepackage[applemac]{inputenc}  
\usepackage[frenchb]{babel}
\usepackage{amsfonts,amsmath}
\usepackage{latexsym}
\usepackage{amssymb}
\usepackage{graphicx}
\usepackage{a4wide}
\usepackage{textcomp}

\def\R{\mathbb R}  
\def\C{\mathbb C}  
\def\H{\mathbb H} 

\newtheorem{theo}{Theorem}

\newtheorem{prop}{Proposition}
\newtheorem{defi}{Definition}
\newtheorem{cor}{Corollary}

\begin{document}
\begin{center}
\textbf{Totally reducible holonomies of torsion-free affine connections}\\
\textbf{Lionel B\'ERARD BERGERY}\\
\end{center}
\begin{center}
En hommage à Marcel BERGER\\
\end{center}
{\bf Summary}

That announcement gives the structure of {\bf totally reducible linear Lie algebras} which are
the {\bf Lie algebra of the holonomy group} of (at least) one torsion-free connection. The result uses the 
(already known) classification of the irreducible ones and some previous (unpublished) works by the author
giving the classification for the pseudo-riemannian totally reducible case.

One describes those Lie subalgebras through a general structure theorem involving
 two constructions and some lists. These constructions give new examples of non irreducible
totally reducible holonomy algebras and also recover some irreducible ones which seem missing
in the previous classification. 

\bigskip
\noindent{\bf 1. Introduction}

A torsion-free linear connection (on a connected n-dimensional manifold $M$) defines 
a holonomy group $H \subset Gl(n,\R)$ (up to conjugacy class). The Lie subalgebra $h\subset gl(n,\R)$ of
$H$ is called the holonomy algebra. It is more precisely the Lie subalgebra of the restricted
holonomy group $H_0$ and it characterizes that connected Lie subgroup of $Gl(n,\R)$. 
For more informations on the holonomy groups, see for example the chapter 10 in the book [Bes]. 

Throughout the paper, for short, a linear Lie algebra which is the holonomy algebra of at least one torsion
free connection will be called a "holonomy". Not any linear Lie algebra is a "holonomy", 
so it is a natural question to investigate which linear Lie algebras are indeed holonomies.
The answer is "known" in the case where $h$ is an irreducible subalgebra i.e., if the linear representation
of $h$ in $\R^n$ (given by the inclusion $h \subset gl(n,\R)$) is irreducible. The "final" classification 
was given by S. Merkulov and L. Schwachhöfer [M-S], after the works of many authors, starting with the
seminal paper by M. Berger [Be1]. A detailed history of the intermediate contributions is explained in two
survey papers by R. Bryant [Br1,2].

Here a "totally reducible" linear Lie algebra is a Lie subalgebra $h$ of a linear algebra $gl(n,\R)$ such
that the linear representation of $h$ on $\R^n$ is totally reducible (or equivalently semi-simple) i.e., is
the direct sum of irreducible representations. Notice that the Lie algebra $h$ is reductive, but not necessarily 
semi-simple. 

Now many holonomies are obtained only for linear connections of locally symmetric spaces.
M. Berger already gave in [Be2] a complete classification of ($1$-connected) symmetric spaces with
totally reducible holonomy. Some of them have a (totally reducible) non irreducible holonomy and are
not products of symmetric spaces. Now one will focus on the non-symmetric case and a "non-symmetric" 
holonomy will be the holonomy of at least one non locally symmetric torsion free connection.\\
{\bf Remarks: 1)} One will {\bf not} follow the usual presentation of the lists which were given in 
the irreducible case. Indeed our constructions give a different insight on those lists.\\
{\bf 2)} The associated "structures" (for example pseudo-riemannian or conformal) will not be studied 
here in order to shorten that announcement.\\ 
{\bf 3)} Notice also that the proof of the results below uses the previous classifications of the irreducible
case and that paper is not a new proof in that case. On the other hand, there are some irreducible
holonomies which were missed in the previous classifications (see the remark at the end of paragraph 8).\\
{\bf 4)} There are some holonomies which are the holonomy algebra of both
symmetric spaces and non locally symmetric manifolds. There are quite few of them in the totally
reducible case and they are all pseudo-riemannian (or "metric") holonomies.

Paragraph 2 contains some definitions and conventions in order to avoid any misleadings, and
paragraph 3 contains the structure theorem.

\bigskip
\noindent{\bf 2. Some definitions for totally reducible linear algebras}
\begin{defi} {\bf : Indecomposable linear Lie algebras}\\
{\bf 1)} Let $h_1 \subset gl(n_1,\R)$ and $h_2 \subset gl(n_2,\R)$ two linear Lie subalgebras.
The direct product of $h_1$ and $h_2$ is the Lie algebra $h = h_1 \oplus h_2$, included 
as a subalgebra of $gl(n,\R)$ with $n = n_1 + n_2$ through the following formula : 
for all $ A \in h_1, B \in h_2, x \in \R^{n_1}, y \in \R^{n_2}$,\\
we have $(A,B)(x,y) = (Ax,By) \in \R^n = \R^{n_1} \times \R^{n_2}.$\\
{\bf 2)} A linear Lie subalgebra is {\bf indecomposable} if it is not the direct product of two
linear Lie subalgebras (with $n_1 n_2 \neq 0$).
\end{defi}
\begin{defi} {\bf : Totally complex, real-complex and totally real linear Lie algebras}\\
{\bf 1)} A linear algebra $h$ is called {\bf totally complex} if it is a {\bf complex} Lie subalgebra of the 
complex Lie algebra $gl(m,\C)$.\\
{\bf 2)} A linear algebra $h$ is called {\bf real-complex} if it is a real {\bf non-complex} Lie subalgebra
of the complex Lie algebra $gl(m,\C)$.\\
{\bf 3)} A linear algebra $h$ is called {\bf totally real} if it is a Lie subalgebra of $gl(n,\R)$ 
and, if $n = 2m$, it is not a subalgebra of the complex Lie algebra $gl(m,\C)$ (even up to conjugacy).
\end{defi}
{\bf Remarks: 1)} The Lie algebra $gl(m,\C)$ is a (real) Lie subalgebra of $gl(2m,\R)$, hence 
a totally complex linear Lie algebra (in  complex dimension $m$) is also a (real) linear Lie algebra 
(in real dimension $n = 2m$).\\
{\bf 2)} The above denominations are used in that paper because quite often, what is called a
"complex" holonomy may be either a complex or a non complex Lie subalgebra of $gl(n,\C)$.
It is important to distinguish carefully theses two cases since the proof of the main theorem use
the "total complexification" for linear Lie algebras, and that may introduce some confusions with the
usual complexification in representation theory.\\
In that direction, one recalls the following well known:
\begin{prop}{\bf : Complexification of indecomposable linear Lie algebras}\\
Let $h \subset gl(n,\R)$ be an indecomposable linear Lie algebra. Then its {\bf total complexification}
$h^{\C} = (h \otimes \C)  \subset gl(n,\C)$ is indecomposable if and only if $h \subset gl(n,\R)$
is {\bf not} a totally complex linear Lie algebra.
\end{prop}

\vfill\eject
\noindent{\bf 3. Structure of totally reducible holonomies}

With all the above preliminaries, one may state the following structure theorem: 
\begin{theo} {\bf : Structure of totally reducible holonomies}\\
{\bf (1) Decomposition:} Any totally reducible holonomy in dimension $n \geqslant 2$ is a 
direct product of linear Lie algebras which are :\\
- either {\bf indecomposable} totally reducible {\bf holonomies},\\
- or $\R = gl(1,\R)$ in dimension $1$ (which is not a holonomy).\\
{\bf (2) Indecomposable case:} Any {\bf totally reducible, indecomposable, non symmetric} holonomy
$h$ is isomorphic (up to conjugacy) to a linear Lie algebra in one of the five "types" I, II, III, IV, V 
described below, where type I, III, IV, V holonomies are given by lists I, III, IV, V in paragraphs
4, 8, 9,10 respectively, and type II holonomies are given by one of the two constructions in
paragraphs 5 and 6.
\end{theo}

\noindent{\bf Types of indecomposable totally reductive holonomies}\\
For the proof of the theorem, it is easier to study together the holonomies with the same properties.
Here are the "characteristics" of these 5 "types":\\
{\bf Type I holonomies} are irreducible linear Lie algebras with non-zero first prolongation and some
further properties. They are the {\bf list I} in paragraph 4. They form the "building blocks" for the
constructions of type II holonomies in paragraph 5.\\
{\bf Type II holonomies} are indecomposable (totally reducible) holonomies which are given by {\bf two
constructions} described in the propositions of paragraphs 5 and 6. As a result, the decomposition as
a direct sum of irreducible factors for an indecomposable holonomy representation may have any
number of irreducible factors. But there are some irreducible type II holonomies, which were already
known. One will not give their list, since it may be deduced easily from the list I (see paragraph 7).\\
{\bf Type III holonomies} are indecomposable holonomies with a special property for the
representation of the semi-simple factor in $h$. They are the {\bf list III} in paragraph 8.\\
{\bf Type IV holonomies} are irreducible symplectic holonomies which are not in previous
lists, and their classification is related with that of complex simple Lie algebras.
They are the {\bf list IV} in paragraph 9.\\
{\bf Type V holonomies} are 10 exceptional irreducible holonomies in low dimensions which
will complete the lists of indecomposable holonomies in the structure theorem.
They are the {\bf list V} in paragraph 10.\\

\noindent{\bf Remarks: 1)} In part (1) of the theorem the condition $n \geqslant 2$ is necessary since 
$\R = gl(1,\R)$ is not a holonomy. Indeed, in {\bf real} dimension $1$ the holonomy is always the
"null" algebra $h = \{ 0 \} \subset \R = gl(1,\R)$. Here all connections are locally symmetric and
$\{ 0 \} \subset gl(1,\R)$ is a irreducible (symmetric) holonomy.\\
 {\bf 2)} If $h$ is the totally reducible holonomy algebra of a {\bf locally symmetric} space, one knows that
it is a direct product of indecomposable totally reducible holonomy algebras of symmetric spaces.
And $\R = gl(1,\R)$ cannot be one of the factors, but $\{0\} \subset gl(1,\R)$ may be a factor.\\
{\bf 3)} Conversely, if $h$ is a {\bf non symmetric} totally reducible holonomy algebra, necessarily the real
dimension $n$ is $\geqslant 2$. And then at least one of its indecomposable factors is non symmetric
or is the exception $\R = gl(1,\R)$.\\
{\bf 4)} Notice that the holonomies in the lists in paragraphs 4, 9, 10 are irreducible. On the other hand, 
the two constructions in paragraphs 5 and 6 and the lists III provide both irreducible
and non irreducible indecomposable holonomies.\\
{\bf 5)} Notice the following consequence of the structure theorem:
\begin{cor} {\bf : Complexification of totally reducible holonomies} \\
A totally reducible linear Lie algebra $h$ {\bf in dimension} $\mathbf{n \geqslant 2}$ is a
holonomy if and only if its (total) complexification $h^{\C} = (h \otimes \C)  \subset gl(n,\C)$
is a totally complex totally reducible holonomy.
\end{cor}

\bigskip
\noindent{\bf 4. Type I holonomies : Irreducible holonomies with $\mathbf{h^{(1)} \neq 0}$ 
and property $\mathcal{C}$.}

Here $h^{(1)} = ((\R^n \circ \R^n)^{\ast} \otimes \R^n) \cap ((\R^n)^{\ast} \otimes h) \subset
((\R^n \otimes \R^n)^{\ast} \otimes \R^n)$ is the first prolongation of the linear Lie algebra $h$.
The irreducible linear Lie algebras $h \subset gl(n,\R)$ {\bf with} $\mathbf{h^{(1)} \neq 0}$ were
classified by S. Kobayashi and T. Nagano [K-N], after a previous work by E. Cartan [Ca].
For future use in the constructions of paragraphs 5 and 6, one selects among them those which 
satisfy a further property $\mathcal{C}$, defined below :
\begin{defi} {\bf : Property} $\mathcal{C}$\\
Let $h \subset gl(n,\R)$ be a linear Lie subalgebra and $h^{(1)}$ its first prolongation.
Denote by $C(h)$ the subspace of $h$ generated by all the maps $B(x, . ) \in h$ for 
all $x \in \R^n$ and all $B \in h^{(1)}$ (viewed as bilinear maps).
Then $h \subset gl(n,\R)$ {\bf satisfies property} $\mathcal{C}$ if  $C(h) = h$ and $h$ has
a non trivial center. 
\end{defi} 
It is an easy exercise to check which ones of the irreducible linear Lie algebras satisfy the above
property $\mathcal{C}$. They are exactly the linear Lie algebras of the list I (A and B) below and
the "exception" $\R = gl(1,\R)$. With property $\mathcal{C}$, the (non trivial) center of $h$ has to
be $\C$ in the totally complex case and $\R$ in the real case.
The following result is already well-known :
\begin{prop}{\bf : Type I holonomies (irreducible holonomies with property $\mathcal{C}$)}\\
(i) All the linear Lie algebras in list {\bf I-A and I-B} below are irreducible, with non-zero first
prolongation, satisfy property $\mathcal{C}$ and they are {\bf non symmetric holonomies}.\\
(ii) The linear algebra $\R = gl(1,\R)$ in dimension $1$ is irreducible, with non-zero first
prolongation, satisfies property $\mathcal{C}$ but it is {\bf not a holonomy}.
\end{prop}
{\bf Remark:} All conventions and notations for the lists below are gathered in paragraph 11 at
the end of the paper.

\noindent{\bf List I-A : Totally complex type I holonomies}
$$\begin{array}{|c||c|c|c|c|}
	\hline
	 &\mathrm{complex}\, h & \mathrm{complex}\, \rho & dim^{\C} & \mathrm{conditions} \\
	\hline
	\hline
	(1) & gl(m,\C) & can & m & m \geqslant 1 \\
	\hline
	(2) & \C \oplus so(m,\C) & \gamma \otimes_{\C} can & m & m \geqslant 3 \\
	\hline
	(3) & \C \oplus sl(p,\C) \oplus sl(q,\C) & \gamma \otimes_{\C} can_p \otimes_{\C} can_q 
	 &  p q & p \geqslant q \geqslant 2 \\
	 & & & & (p,q) \neq (2,2) \\
	\hline
	(4) & gl(m,\C) & Sym^2(can) & \frac{m(m+1)}{2} & n \geqslant 3 \\
	\hline
	(5) & gl(m,\C) & Ext^2(can) & \frac{m(m-1)}{2} & n \geqslant 5 \\
	\hline
	(6) & \C \oplus so(10,\C) & \gamma \otimes_{\C} (half \, spin) & 16 &\\
	\hline
	(7) & \C \oplus E_6^{\C} & \gamma \otimes_{\C} can & 27 & \\
	\hline
\end{array}$$
\noindent{\bf List I-B : Totally real type I holonomies}
$$\begin{array}{|c||c|c|c|c|}
	\hline
	 & \mathrm{real}\, h & \mathrm{totally\, real}\, \rho & dim^{\R} & \mathrm{conditions} \\
	\hline
	\hline
	(1a) & gl(n,\R) & can & n & n \geqslant 2 \\
	\hline
	(2a) & \R \oplus so(p,q) & \gamma \otimes can & p+q & p \geqslant q \geqslant 0 \\
	 & & & & p+q \geqslant 3 \\
	\hline
	(3a) & \R \oplus sl(p,\R) \oplus sl(q,\R) & \gamma \otimes can_p \otimes can_q 
	 & p q & p \geqslant q \geqslant 2 \\
	 & & & & (p,q) \neq (2,2) \\
	\hline
	(3b) & \R \oplus sl(p,\C) & \gamma \otimes Herm(can_p) &  p^2 & p \geqslant 3 \\
	\hline
	(3c) & \R \oplus sl(p,\H) \oplus sl(q,\H) & \gamma \otimes (can_p \otimes_{\H} can_q)
	 & 4 p q & p \geqslant q \geqslant 1 \\
	  & & & & (p,q) \neq (1,1) \\
	\hline
	(4a) & gl(m,\R) & Sym^2(can) & \frac{m(m+1)}{2} & m \geqslant 3 \\
	\hline
	(4b) & gl(m,\H) & Antiherm(can) & m(2m+1) & m \geqslant 2 \\
	\hline
	(5a) & gl(m,\R) & Ext^2(can) &  \frac{m(m-1)}{2} & m \geqslant 5 \\
	\hline
	(5b) & gl(m,\H) & Herm(can) & m(2m-1) & m \geqslant 3 \\
	\hline
	(6a) & \R \oplus so(5,5) & \gamma \otimes (half \, spin^{\R}) & 16 & \\
	\hline
	(6b) & \R \oplus so(9,1) & \gamma \otimes (half \, spin^{\R}) & 16 & \\
	\hline
	(7a) & \R \oplus E_6^1 & \gamma \otimes can & 27 & \\
	\hline
	(7b) & \R \oplus E_6^4 & \gamma \otimes  can & 27 & \\	
	\hline
\end{array}$$

\bigskip
\noindent{\bf 5. Type II holonomies : totally complex construction}

Now two constructions will give a lot of other indecomposable totally reducible holonomies
starting from the list I in paragraph 4 (and also $\R$). Most of them are not irreducible, but some are.
Notice that in the non-irreducible case, the center of $h$ gives the indecomposable property for the
holonomy. First, here is the complex construction which gives totally complex holonomies.\\

\noindent{\bf Data for the complex construction}

Let  $h_j \subset gl(n_j,\C),\, j = 1, ... , p$ be $p$ complex linear Lie algebras from the
list 1-A above ($p \geqslant 1$). Such a Lie algebra may be written $h_j = a_j \oplus s_j$,
where $s_j$ is semi-simple and $a_j$ is the center. More precisely, using the linear
representation, $a_j = \C Id_j$, where $Id_j \in gl(n_j,\C)$ is the identity map.
Denote by $a$ the complex $p$-dimensional abelian Lie algebra which is the direct
product of all those $1$-dimensional centers i.e., $a = \oplus_{j = 1}^{j = p} a_j$.\\
The basis of the complex $m$-dimensional vector space $a$ given by the $Id_j$
identifies $a$ with $\C^p$. Inside $a$, define a "generic" complex hyperplane $z$
by one homogeneous equation\\
$z = \{(w_1, ..., w_p) \in \C^p = a \, ; \, \sum_{j = 1}^{j = p} \lambda_j w_j = 0 \},$
with the following property $\mathbf{(\mathcal{P}_1)}$: 

\noindent{\bf Property} $\mathbf{(\mathcal{P}_1)}$: {\bf all the complex numbers} 
$\mathbf{\lambda_j}\, (j = 1, ... , p)$ {\bf are non zero}.\\
(Obviously, for any given complex hyperplane $z$ 
the complex numbers $\lambda_j$ are unique up to a non zero common factor.)

Let $h$ be the complex Lie algebra which is the direct product of the complex abelian
Lie algebra $z$ and all the semi-simple $s_j$ above i.e., $h = z \oplus s$ where
$s = \oplus _{j = 1}^{j = p} s_j$.

Let $m = n_1 + ... + n_p$. The vector space $\C^m$ may be viewed as 
the direct product $\oplus_{j = 1}^{j = p} \C^{n_j}$ of all the representation spaces of
the $h_j$. Then $h$ becomes a complex Lie subalgebra of $gl(m,\C)$ through
the representation $\rho$ given by the following formula: \\
for all $ \, (w_1, ..., w_p) \in z \, , \, (A_1, ..., A_p) \in s \, ,
 \, (X_1, ..., X_p) \in \C^m ,$\\
we have $\rho(w_1, ..., w_p,A_1, ..., A_p)(X_1, ..., X_p) =
(w_1 \, X_1 + A_1 (X_1), ... ,w_p \, X_p + A_p (X_p)).$
\begin{prop} {\bf : "Totally complex Type II holonomies"}\\
With the above "data for the complex construction", $h \subset gl(m,\C)$ is a {\bf totally complex
type II indecomposable totally reducible holonomy}.
\end{prop}
\noindent{\bf Remark:} In that construction, the linear representation $\rho$ of the Lie algebra $h$
in the vector space $\C^m$ is the direct {\bf sum} of $p$ irreducible representations. But it is not 
a direct {\bf product} since $z$ is "generic" in $a$.

\bigskip
\noindent{\bf 6. Type II holonomies : real construction}

A similar construction applies to the real case and gives the totally real forms
of the totally complex linear algebras given by the complex construction.
Once again the real $1$-dimensional case plays a special role.\\

\noindent{\bf Data for the real construction}

Let  $h_j \subset gl(n_j,\C),\, j = 1, ..., p$ be $p$ {\bf complex} linear Lie algebras from the
list 1-A above. Let $h_k \subset gl(n_k,\R),\, k = p+1, ..., p+q$ be $q$ (real) linear Lie algebras
which are either from the list 1-B above or $\R = gl(1,\R)$.
Such a Lie algebra may be written $h_i = a_i \oplus s_i$, where
$s_i$ is semi-simple and $a_i $ is the center.
More precisely, $a_i = \C Id_i$ or $\R Id_i$, where $Id_i$ is the identity in $gl(n_i,\C)$
or $gl(n_i,\R)$ respectively.

Let $p \geqslant 0,\, q \geqslant 0$ and $m = 2p+q \geqslant 1$.
Denote by $a$ the {\bf real} $m$-dimensional abelian Lie algebra which is the direct product
of all the centers i.e., $a = \oplus_{i = 1}^{i = p+q} a_i$.
There is a basis for $a$ given by the $Id_j$ and $i\,  Id_j$ for $j = 1,...,p$ and the $Id_k$ for 
$k = p+1, ..., p+q$. That basis identify the real $m$-dimensional vector space $a$ with 
$\R^m = \C^p \oplus \R^q$.

Inside $a$, let $z$ be a {\bf real hyperplane} given by one homogeneous equation\\
$z = \{(w_1, ..., w_p,u_{p+1},...,u_{p+q}) \in  a = \C^p  \times \R^q = \R^m \, ; \, 
\sum_{j = 1}^{j = p} Im(\lambda_j w_j) + \sum_{k = p+1}^{k = p+q} \mu_{k} u_k)= 0 \},$\\
with the following property $\mathbf{(\mathcal{P}_2)}$:\\
{\bf Property} $\mathbf{(\mathcal{P}_2)}$: 
{\bf all the complex numbers} $\mathbf{\lambda_j}\, (j = 1, ... , p)$
{\bf and all the real numbers} $\mathbf{\mu_k}\, (k = p+1, ... , p+q)$
{\bf are non zero}.\\
(Here $z$ is a "generic" real hyperplane in $a$. Obviously, for any given hyperplane $z$ 
the complex numbers $\lambda_j$ and the real numbers $\mu_k$ are unique up to a non
zero {\bf real} factor.)

Let $h$ be the (non complex) Lie algebra which is the direct product of the abelian
Lie algebra $z$ and all the semi-simple $s_i$ above i.e., $h = z \oplus_{i = 1}^{i = p+q} s_i$.\\
Let $n = 2(n_1+ ... + n_p) + (n_{p+1} + ... + n_{p+q})$.
Assume $n \geqslant 2$ in order to avoid the $1$ dimensional case.
The vector space $\R^n$ may be viewed as the direct product 
$\oplus_{j = 1}^{j = m} \C^{n_j} \oplus \oplus_{k = p+1}^{k = p+q} \R^{n_k}$ of all the
representation spaces of the $h_i$'s. 

Then $h$ becomes a real Lie subalgebra of $gl(n,\R)$ through
the representation $\rho$ given by the following formula:\\
for all $(w_1, ..., w_p,u_{p+1},...,u_{p+q}) \in z \, , \, (A_1, ..., A_{p+q})
\in s \,, \, (X_1, ... , X_p, Y_{p+1},...,Y_{p+q}) \in \R^n$,\\
we have 
$\rho(w_1, ..., w_p,u_{p+1},...,u_{p+q},A_1, ..., A_{p+q}) (X_1, ... , X_p,Y_{p+1},...,Y_{p+q}) =\\
(w_1 \, X_1 + A_1\, (X_1), ... ,w_p \, X_p + A_p (X_p),
u_{p+1}\, Y_{p+1} + A_{p+1} (Y_{p+1}), ... , u_{p+q}\, Y_{p+q} + A_{p+q} (Y_{p+q}))$.
\begin{prop} {\bf : "Real (non complex) Type III holonomies"}\\
With the above "data for the real construction", $h \subset gl(n,\R)$ is a real (non complex)
{\bf type II indecomposable totally reducible holonomy}.
\end{prop}
\noindent{\bf Remarks:} In that construction, the linear representation $\rho$ of the Lie algebra
$h$ in the vector space $\R^n$ is the direct {\bf sum} of $p+q$ irreducible representations.
But it is not a direct {\bf product} since $z$ is "generic" in $a$.
When $q = 0$ then the Lie algebra $h$ is real since $z$ is a real hyperplane, but the holonomy
representation $\rho$ is complex. Then we have a real-complex holonomy.
On the other hand, when $p = 0$ all the irreducible factors are totally real.

\bigskip
\noindent{\bf 7. Two remarks on type II holonomies}

\smallskip
\noindent{\bf (1) Irreducible type II holonomies}

Notice that applying the above constructions {\bf with only one factor} give all the {\bf irreducible
type II holonomies}. First, here are some examples:
$$\begin{array}{|c||c|c|c|c|c|}
	\hline
	(A) & sl(m,\C) & can & dim^{\C} = m & m \geqslant 2 & \mathrm{totally\, complex}\\
	\hline
	(B) & sl(n,\R) & can & dim^{\R} = n & n \geqslant 2 & \mathrm{totally\, real} \\
	\hline
	(C) & \R \oplus sl(m,\C) & \sigma_{\theta} \otimes_{\C} can & dim^{\C} = m & m \geqslant 1
	 &  \mathrm{real-complex}\\	
	\hline
\end{array}$$
Now the general construction of type II {\bf irreducible} holonomies is the following :\\
{\bf II-A Totally complex irreducible type II holonomies:} they result from the complex
construction with only one factor and there is one such example $h \subset gl(n,\C)$ for each
member $\C \oplus h \subset gl(n,\C)$ of List I-A. (Do not forget $gl(,n\C) = \C \oplus sl(n,\C)$).\\
{\bf II-B Totally real irreducible type II holonomies:} they result from the real construction
with only one real non complex factor and there is one such example $h \subset gl(n,\R)$ for each
member $\R \oplus h \subset gl(n,\R)$ of List I-B.\\
{\bf II-C Real-complex irreducible type II holonomies:} they result from the real construction
with only one totally complex factor and there is one such example $\R \oplus h \subset gl(n,\C)$
for each member $\C \oplus h \subset gl(n,\C)$ of List I-A and representation $\sigma_{\theta}$
for the center $\R$ with $ 0 \leqslant \theta \leqslant \frac{\pi}{2}$.\\
In those ways, one gets all the irreducible type II holonomies.\\
\noindent{\bf Remark:} Notice that in the classical lists for the irreducible case in [Br1], [Br2] or
[M-S], they are listed with the same semi-simple factor of $h$ and various factors for the center.\\

\noindent{\bf (2)} The following corollary is an obvious consequence of the structure theorem :
\begin{cor} {\bf : Totally reducible indecomposable holonomies with many factors}\\
If a totally reducible indecomposable holonomy has at least 3 irreducible factors in the decomposition
of the holonomy in a direct sum, then it is a type II holonomy.
\end{cor}
Notice that there are other (non type II) indecomposable totally reducible holonomies whose
decomposition in direct sum contains precisely two factors. They are type III  holonomies below
in paragraph 8 (list III A and B).

\bigskip
\noindent{\bf 8. Type III holonomies : indecomposable holonomies with property $\mathcal{S}$}

As in that title, type III holonomies are indecomposable totally reducible holonomies such that
the restriction of the holonomy representation to the semi-simple factor in $h$ satisfies the
following (called property $\mathcal{S}$ for short): 
\begin{defi} {\bf : Property} $\mathcal{S}$\\
Let $h \subset gl(n,\R)$ be an indecomposable totally reducible holonomy and let $s$ be the
semi-simple factor in $h$. Denote by $\rho^{s s}$ the restriction of the holonomy representation
$\rho$ to $s$. \\
Then $h$ satisfies property $\mathcal{S}$ if and only it is not of type I or II and it satisfies one of
the following:\\
(i) $\rho^{s s}$ is not irreducible,\\
(ii) $s$ is non complex and $\rho^{s s}$ is complex,\\
(iii) $s = s_1 \oplus s_2$, with $s_2 = sl(2,\C)$ or $su(2)$ or $sl(2,\R)$ and $\rho^{s s}$ is the 
tensor product of a representation of $s_1$ with the canonical representation of $s_2$.
\end{defi} 
{\bf Remark: } The condition (iii) is some sort of "global" quaternionic or paraquaternionic 
structure for the holonomy representation (cases (7) of List III).

\medskip
In the totally complex case one get only (i) or (iii) in property $\mathcal{S}$. 
In the case (iii), the representation $\rho$ is irreducible. In the case (i),
the direct sum decomposition of $\rho^{s s}$ has exactly 2 factors i.e., 
$\rho^{s s} = \rho^{s s}_1 \oplus \rho^{s s}_2$.
In the lists below, there are different properties for these two restricted representations\\
- in the cases (1), (2), (3) and (6):  $\rho^{s s}_2 = \rho^{s s}_1$,\\
- in the cases (4), (5), (6) and also, if m = 2, in the cases (1), (2), (3): 
$\rho^{s s}_2 = (\rho^{s s}_1)^{\ast} $\\
- in the cases (8), (9), (10), (11), $\rho^{s s}_1$ and $\rho^{s s}_2$ have different dimensions.\\

\noindent{\bf List III-A : Totally complex type III holonomies}
$$\begin{array}{|c||c|c|c|c|}
	\hline
	 & \mathrm{complex}\, h & \mathrm{complex}\, \rho & dim^{\C}
	 & \mathrm{conditions} \\
	\hline
	\hline
	(1) & sl(m,\C) &  can \oplus can & 2m & m \geqslant 2 \\	   
	\hline	                 	
	(2) & \C \oplus sl(m,\C) & (\gamma^{\C} \otimes_{\C} can) \oplus (\gamma^{\C}(\lambda)
	\otimes_{\C} can) & 2m & m \geqslant 2 \\	
	\hline
	(3) & \C^2 \oplus sl(m,\C) & (\pi_1^{\C} \otimes_{\C} can) \oplus (\pi_2^{\C} \otimes_{\C} can)
	 & 2m & m \geqslant 2 \\
	\hline
	\hline
	(4)  & sl(m,\C) &  can \oplus can^{\ast} & 2m & m \geqslant 3 \\	
	\hline
	(5) & gl(m,\C) &  can \oplus can^{\ast} & 2m & m \geqslant 3 \\
	\hline
	\hline
	(6)  & sp(m,\C) &  can \oplus can & 4m & m \geqslant 2 \\	
	\hline
	(7) & sp(m,\C) \oplus sp(1,\C) & can \otimes_{\C} can  & 4m & m \geqslant 2 \\	
	\hline
	\hline
	(8) & \C \oplus sl(p,\C) & (\gamma^{\C} \otimes_{\C} \rho_0^{\C} \otimes_{\C} can_q)\,
	\oplus & (p+1)q  & p \geqslant 2 \\
	 &  \oplus \, sl(q,\C) & (\gamma^{\C}(\frac{p}{p+q}) \otimes_{\C} can_p \otimes_{\C} can_q)
	  & & q \geqslant 2 \\
	\hline
	(9) & \C^2 \oplus sl(m,\C) & (\pi_1^{\C} \otimes_{\C} \rho_0^{\C} \otimes_{\C} can_2) \, \oplus 
	 & 2m + 2 & m \geqslant 2 \\
	 & \oplus \, sl(2,\C)  & (\pi_2^{\C} \otimes_{\C} can_m \otimes_{\C} can_2) & & \\
	\hline
	\hline
	(10) & \C \oplus sl(m,\C) & (\gamma^{\C} \otimes_{\C}  can)\, \oplus & \frac{1}{2}m(m+3) 
	& m \geqslant 2  \\
	 & & (\gamma^{\C}(\frac{1}{2}) \otimes_{\C} Sym^2(can)) & & \\
	\hline
	(11) & \C^2 \oplus sl(2,\C) & (\pi_1^{\C} \otimes_{\C} can)\, \oplus & 5 &\\
	 & & (\pi_2^{\C} \otimes_{\C} Sym^2(can)) & & \\
	\hline
\end{array}$$
\noindent{\bf Liste III-B : Totally real type III holonomies}
$$\begin{array}{|c||c|c|c|c|}
	\hline
	 & \mathrm{real}\, h & \mathrm{totally\, real}\, \rho
	  & dim^{\R} & \mathrm{conditions} \\
	\hline
	\hline
	(2a) & \R \oplus sl(m,\R) & (\gamma^{\R} \otimes can) \oplus
	(\gamma^{\R}(\mu) \otimes can) & 2m & m \geqslant 2 \\	
	 & & & & \mu \neq 1 \\
	\hline
	(3a) & \R^2 \oplus sl(m,\R) &  (\pi_1^{\R} \otimes can) \oplus (\pi_2^{\R} \otimes can) 
	 & 2m & m \geqslant 2 \\
	\hline
	\hline
	(4a) & sl(m,\R) &  can \oplus can^{\ast} & 2m & m \geqslant 3 \\	
	\hline
	(5a) & gl(m,\R) &  can \oplus can^{\ast} & 2m & m \geqslant 3 \\
	\hline
	\hline
	(7a) & sp(m,\R) \oplus sp(1,\R) & can \otimes can  & 4m & m \geqslant 2 \\	
	\hline
	(7b) & sp(p,q) \oplus sp(1) & can \otimes_{\H} can  & 4(p+q) & p+q \geqslant 2 \\	
	 & & & & p \geqslant q \geqslant 0 \\
	\hline
	\hline
	(8a) & \R \oplus sl(p,\R)& (\gamma^{\R} \otimes \rho_0 \otimes can_q)\,
	\oplus & (p+1)q  & p \geqslant 2 \\
	 &  \oplus\, sl(q,\R)  & ( \gamma^{\R}(\frac{p}{p+q}) \otimes can_p \otimes_{\C} can_q)
	& & q \geqslant 2 \\
	\hline
	(9a) & \R^2 \oplus sl(m,\R) & (\pi_1^{\R} \otimes \rho_0 \otimes can_2)\, \oplus \
	& 2m + 2 & m \geqslant 2 \\
	 & \oplus \, sl(2,\R) & (\pi_2^{\R} \otimes can_m \otimes can_2)  & & \\
	\hline
	\hline
	(10a) & \R \oplus sl(m,\R) & (\gamma^{\R} \otimes can) \, \oplus & \frac{1}{2}m(m+3) 
	& m \geqslant 2 \\
	 & & (\gamma^{\R}(\frac{1}{2})  \otimes Sym^2(can)) & & \\
	\hline
	(11a) & \R^2 \oplus sl(2,\R) & (\pi_1^{\R} \otimes can)\, \oplus & 5 & \\	 
	 & & (\pi_2^{\R} \otimes Sym^2(can)) & & \\
	\hline
\end{array}$$
\noindent{\bf Liste III-C : Real-complex type III holonomies}
$$\begin{array}{|c||c|c|c|c|}
	\hline
	 & \mathrm{real}\, h & \mathrm{complex}\, \rho & dim^{\C} & \mathrm{conditions} \\
	\hline
	\hline
	(1a) & sl(m,\R) &  can \oplus can & m & m \geqslant 2 \\
	 & & = can \otimes \C & & \\	
	\hline
	(1b) & sl(m,\H) & can & 2m & m \geqslant 1 \\
	\hline
	(2b) & \R \oplus sl(m,\R) &   \sigma_{\theta}^{\C} \otimes can & m & m \geqslant 2 \\
	 & & & &  0 \leqslant \theta \leqslant \frac{\pi}{2} \\
	\hline
	(2c) & \R \oplus sl(m,\H) &  \sigma_{\theta}^{\C} \otimes_{\C} can & 2m
	 & m \geqslant 1 \\
	 & & & &  0 \leqslant \theta \leqslant \frac{\pi}{2} \\
	\hline
	(3b) & \C \oplus sl(m,\R) & \gamma^{\C} \otimes can & m & m \geqslant 2\\	
	\hline
	(3c) & \C \oplus sl(m,\H) & \gamma^{\C} \otimes_{\C} can & 2m & m \geqslant 1 \\
	\hline
	\hline
	(4b) & su(p,q) & can & (p+q) & p \geqslant q \geqslant 0 \\
	 & & & & p+q \geqslant 3 \\
	\hline
	(5b) & u(p,q) & can & (p+q) & p \geqslant q \geqslant 0 \\ 
	 & & & & p+q \geqslant 3 \\
	\hline
	\hline
	(6a) & sp(m,\R) &  can \oplus can & 2 m & m \geqslant 2 \\	
	 & & = can \otimes \C & & \\	
	\hline
	(6b) & sp(p,q) & can & 2(p+q) & p \geqslant q \geqslant 0 \\
	 & & & & p+q \geqslant 2 \\
	\hline
\end{array}$$
{\bf Remarks: 1)} The cases (2b) (with $\theta \neq 0$), (2c), (3b), (3c) give irreducible
holonomies. They appear in the previous classifications (for the irreducible case)
only for $m = 2$ in the cases (2b) and (3b) and for $m = 1$ in the cases (2c) and (3c). 
That some of them may exist also for higher $m$ was already noticed by R. Bryant [Br3]
and D. Joyce. But (it seems to me that) all of them are indeed holonomies.\\
{\bf 2)} Notice that the cases (4), (5), (6), (4a), (5a), (6a) (with also the cases $m =2$
of (1), (1a) and the cases $m = 2$ and $\lambda = -1$ of (2) and 
$m = 2$ and $\mu = -1$ of(2a)) are the non irreducible indecomposable totally reducible
pseudo-riemannian holonomies already studied in [BB-I] and [BB].\\

\noindent{\bf 9. Type IV holonomies : special irreducible symplectic holonomies}

In 1996, Q.S. Chi, S.A. Merkulov and L.J. Schwachhöfer discovered new "exotic" irreducible
holonomies [CMS]. They related them with the list of "Wolf spaces" (compact quaternion-Kähler
symmetric spaces), and they studied their properties (with relations to supersymmetry).
It is possible to describe these holonomies in the following way.

Let $g$ be a complex {\bf simple} Lie algebra. A construction due initially to J.A. Wolf [W]
associates to each $g$ a (non simple) complex subalgebra $k = h \oplus sp(1,\C)$ such that
the representation of $h \oplus sp(1,\C)$ on the quotient space $g/k$ is the tensorial product of
a symplectic representation
$\rho$ of $h$ and the canonical 2-dimensional representation of $sp(1,\C)$ (and such a
subalgebra $k$ is unique up to conjugacy in $g$). Then the associated homogeneous space
$G/(H\, Sp(1,\C)$ is symmetric and its "compact real form" is a compact quaternion-kähler
symmetric spaces (they are called  a Wolf's spaces).

Furthermore, with that representation $\rho$, the Lie algebra $h$ becomes a holonomy.
If $g = sl(n+2,\C)$, that construction gives $h = gl(n,\C)$ with $\rho = can \oplus can^{\ast}$,
and such a holonomy is already in List III-A (5). Then the other simple Lie algebras
give the holonomies in the list IV-A below. Notice that $g = sp(n+1,\C)$ gives $h = sp(n,\C)$
with its canonical representation, which was already known to Berger. The other cases
are the new cases discovered by Q.S. Chi, S.A. Merkulov and L.J. Schwachhöfer.

On the other hand, the tensorial product representation for $h \oplus sp(1,\C)$ is usually only
symmetric. There is one family which is both symmetric and non symmetric, with $h = sp(n,\C)$.
This gives the holonomies (7), (7a) and (7b) in the list III.

\bigskip
\noindent{\bf List IV-A :  : Totally complex type IV holonomies}
$$\begin{array}{|c||c|c|c|c|}
	\hline
	 & \mathrm{complex}\, h & \mathrm{complex}\, \rho & dim^{\C} & \mathrm{conditions} \\	
	\hline
	\hline
	(1) & sp(m,\C) & can & 2m & m \geqslant 2 \\
	\hline
	(2) & so(m,\C) \oplus sp(1,\C) & can \otimes_{\C} can & 2m & m \geqslant 3 \\
	\hline
	(3) & sp(1,\C) & Sym^3(can) & 4 &  \\	
	\hline
	(4) & sp(3,\C) & Ext^3_0(can) & 14 &  \\	
	\hline
	(5) & sl(6,\C) & Ext^3(can) & 20 & \\
	\hline
	(6) & so(12,\C) & half \, spin & 32 & \\
	\hline
	(7) & E_7^{\C} & can & 56 & \\ 
	\hline
\end{array}$$
\noindent{\bf List IV-B : Totally real type IV holonomies}
$$\begin{array}{|c||c|c|c|c|}
	\hline
	 & \mathrm{real}\, h & \mathrm{totally\, real}\, \rho & dim^{\R} & \mathrm{conditions} \\
	\hline
	\hline
	(1a) & sp(m,\R) & can & 2m & m \geqslant 2 \\
	\hline
	(2a) & so(p,q) \oplus sp(1,\R) & can \otimes can & 2(p+q) & p+q \geqslant 3 \\
	 & & & & p \geqslant q \geqslant 0 \\
	\hline
	(2b) & so(m,\H) \oplus sp(1) & can \otimes_{\H} can & 4 m & m \geqslant 2 \\
	\hline
	(3a) & sp(1,\R) & Sym^3(can) & 4 & \\	
	\hline
	(4a) & sp(3,\R) &  Ext^3_0(can) & 14 & \\	
	\hline
	(5a) & sl(6,\R) &  Ext^3(can) & 20 & \\
	\hline
	(5b) & su(3,3) &  Ext^3(can)^{\R} & 20 & \\
	\hline
	(5c) & su(5,1) &  Ext^3(can)^{\R} & 20 & \\
	\hline
	(6a) & so(6,6) & half \, spin^{\R} & 32 & \\
	\hline
	(6b) & so(10,2) & half \, spin^{\R} & 32 & \\
	\hline
	(6c) & so(6,\H) & half \, spin^{\R} & 32 & \\
	\hline
	(7a) & E_7^1 & can & 56 & \\ 
	\hline
	(7b) & E_7^3 & can & 56 & \\ 
	\hline
\end{array}$$

\bigskip
\noindent{\bf 10. Type V holonomies : exceptional irreducible holonomies}

In the end, there are 10 irreducible holonomies which do not fit in the above lists or
constructions. They are called exceptional holonomies (obviously a list of "exceptional"
holonomies may depend on various choices for the lists).

Among them, six may be related to the theory of octonions: cases (1), (2) and their real forms.

Finally, there is no "conformally symplectic" holonomy in dimension $n \geqslant 6$,
but such holonomies do exist in dimension $2$ and $4$. Dimension $2$ is already contained
in the previous lists. In dimension $4$, we get the "exceptional conformally symplectic"
holonomies of cases (3), (3a), (4) and (4a).\\

\noindent{\bf List V-A : Totally complex type V holonomies}
$$\begin{array}{|c||c|c|c|c|}
	\hline
	 & \mathrm{complex}\, h & \mathrm{complex}\, \rho & dim^{\C} \\	
	\hline
	\hline
	(1) & G_2^{\C} & can & 7 \\
	\hline
	(2) & so(7,\C) & spin & 8 \\
	\hline
	\hline
	(3) & \C \oplus sl(2,\C) & \gamma^{\C} \otimes_{\C} Sym^3(can) & 4 \\ 	
	\hline
	(4) & \C \oplus sp(2,\C) & \gamma^{\C} \otimes_{\C} can & 4 \\ 	
	\hline
\end{array}$$
\noindent{\bf List V-B : Totally real type V holonomies}
$$\begin{array}{|c||c|c|c|c|}
	\hline
	 & \mathrm{real}\, h & \mathrm{totally\, real}\, \rho & dim^{\R} \\
	\hline
	\hline
	(1a) & G_2^1 & can & 7 \\
	\hline
	(1b) & G_2 & can & 7 \\
	\hline
	(2a) & so(4,3) & spin^{\R} & 8 \\
	\hline
	(2b) & so(7) & spin^{\R} & 8 \\
	\hline
	\hline
	(3a) & \R \oplus sl(2,\R) &  \gamma^{\R} \otimes Sym^3(can) & 4 \\ 	
	\hline
	(4a) & \R \oplus sp(2,\R) &  \gamma^{\R} \otimes can & 4 \\			 
	\hline
\end{array}$$

\bigskip
\noindent{\bf 11. Notations and conventions} 

\noindent {\bf Notations for the lists}\\
In all the lists of that paper, holonomies are described through
the following presentation in 4 columns: (a) $h$, (b) $\rho$, (c) dim, (d) conditions, where\\
(a) $h$ is an "abstract" Lie algebra,\\
(b) $\rho$ is the representation which turns $h$ into a subalgebra of some $gl(n,\R)$ 
(or $gl(m,\C)$),\\
(and the notations and conventions for those $\rho$ are below),\\
(c) $dim^{\R}$ [resp. $dim^{\C}$] is the real [resp. complex] dimension of the
representation space,\\
(d) "conditions" are some conditions on the parameters in order to avoid repetitions\\
(mainly due to classical isomorphisms in low dimensions).

Moreover, all the lists of that paper are divided in 2 or 3 parts (A), (B), (C) according to the
following properties for the holonomies: \\
{\bf (A) totally complex holonomies} i.e., $h$ is a complex Lie subalgebra of $gl(n,\C)$,\\
{\bf (B) totally real holonomies} i.e., the representation space has no $\rho(h)$-invariant
complex structure,\\
{\bf (C) real-complex holonomies} i.e., $h$ is non complex, but the representation $\rho$
is complex.\\

\noindent {\bf Notations for the representation $\mathbf{\rho}$}\\
Usually, the {\bf representation} $\mathbf{\rho}$ is described through a tensorial product of
representations of the center $a$ (if there is one) and each of the simple factors in $h$.
The only exceptions are for the (full) linear algebras $gl(m,\R)$,  $gl(m,\C)$ or  $gl(m,\H)$
which are considered in the same way than simple factors. Also $so(4,\C)$, $so(4)$ and
$so(2,2)$ are sometimes treated in the same way.

If $\rho$ is not irreducible and is the direct sum of two irreducible factors, it is denoted
by $\rho = \rho_1 \otimes \rho_2$ and now each factor is described as above. In the 
special case where $\rho$ is complex-irreducible and not real irreducible, one describes both.

\noindent{\bf Conventions for the center of $h$}\\
The canonical $1$-dimensional representation of $\R = gl(1,\R)$ in $\R$
[resp. of $\C = gl(1,\C)$ on $\C$] is always denoted by $\gamma^{\R}$ [resp. $\gamma^{\C}$].\\
In the case of a $1$-dimensional factor $\C$ in a 2-dimensional (complex) center, the representation
$\gamma^{\C}(\lambda)$ of $\C$ on $\C$ is given by $\gamma^{\C}(\lambda) = \lambda \, Id$,
where the parameter $\lambda$ is a complex number. In order to avoid isomorphisms, one restricts
furthermore to $| \lambda  | \leqslant 1$ and if $| \lambda  | = 1$ then $Im(\lambda) \geqslant 0$.\\
In the case of a $1$-dimensional factor $\R$ in a 2-dimensional (real) center, the representation
$\gamma^{\R}(\mu)$ of $\R$ on $\R$ is given by $\gamma^{\R}(\mu) = \mu \, Id$, where the
parameter $\mu$ is a real number. In order to avoid isomorphisms, one restricts furthermore 
to $| \mu  | \leqslant 1$.\\
The complex representations $\sigma_{\theta}^{\C}$ of $\R$ in the {\bf complex} $1$-dimensional
space $\C$ are given by $\sigma_{\theta}^{\C}(t) = t \, e^{i\,\theta} \, Id$, where  $\theta$ is real
(some angle). In order to avoid isomorphisms, one restricts furthermore to $\theta \in [0, \frac{\pi}{2}]$.
Now $\sigma_{\theta}^{\C}$ is $\R$-irreducible if and only if $\theta \neq 0$.

If  there is a $2$-dimensional center (for a non-irreducible $\rho$), $\pi_1^{\R}$ and $\pi_2^{\R}$
[resp. $\pi_1^{\C}$ and $\pi_2^{\C}$] are the canonical representations of  $\R^2$ on $\R$ 
[resp. of $\C^2$ on $\C$] given by $\pi_1(t,u) = t\, Id$ and $\pi_1(t,u) = u\, Id$.

\noindent{\bf Conventions or the simple factors of $h$}\\
Here the notation $\mathbf{can}$ is used either for any "canonical" representation
of a classical Lie algebra or for the lowest dimensional representation of exotic ones.\\
And $\rho_0^{\C}$ [resp. $\rho_0^{\R}$] is the trivial $1$-dimensional complex [resp. real]
representation i.e., the corresponding simple factor is in the kernel of the representation.\\
Sometimes a complex-irreducible representation $\eta$ of  a real Lie algebra may have a "real
structure" i.e., be the complexification of a real representation, here denoted by $\eta^{\R}$.\\
Notice a very special case: {\bf only one} of the two (complex) half-spin representations of $so(6,\H)$
has a real structure and gives a half-spin$^{\R}$ (and a holonomy).\\
And one recall that $\H^r \otimes_{\H} \H^s = \R^{4 r s}$ is only a real vector space,
with no invariant complex structure for the corresponding (tensorial product) representation.

\bigskip
\noindent{\bf 12. Bibliography}\\

\noindent\textbf{[BB] Lionel B\'ERARD BERGERY}, \textit{L'holonomie des variétés pseudo-riemanniennes:
le cas totalement réductible}, Preprint Institut Elie Cartan (2007) (unpublished).\\
\textbf{[BB-I] Lionel B\'ERARD BERGERY and Aziz IKEMAKHEN}, \textit{Sur l'holonomie des
variétés pseudo-riemanniennes de signature (n,n)},
Bull. Soc. Math. France 125 (1997), 93-114.\\
\textbf{[Be1] Marcel BERGER}, \textit{Sur les groupes d'holonomie des variétés à connexion affine
et des variétés riemanniennes}, Bull. Soc. Math. France 83 (1955), 279-330.\\
\textbf{[Be2] Marcel BERGER}, \textit{Les espaces symétriques non compacts},
Ann. Sci. \'Ecole Norm. Sup. 74 (1957), 85-177.\\
\textbf{[Bes] Arthur L. BESSE}, \textit{Einstein manifolds}, Springer-Verlag (1987).\\
Last edition in "Classics in Mathematics" (2008).\\
\textbf{[Br1] Robert L. BRYANT}, \textit{Classical, exceptional and exotic holonomies: a status report},
Actes de la Table Ronde de Géométrie Différentielle en l'honneur de Marcel Berger,
Sémin. Congr. 1, Soc. Math. de France, (1996), 93-165.\\
\textbf{[Br2] Robert L. BRYANT}, \textit{Recent advances in the theory of holonomy},
Séminaire Bourbaki 51ème année (1998-99) n° 861,
Astérisque 266 n°5 (2000) 351-374.\\
\textbf{[Br3] Robert L. BRYANT}, \textit{Some remarks on Finsler manifolds with constant flag curvature},
Houston Journal of Mathematics, vol. 28 no. 2 (2002), 221-262.\\
\textbf{[Ca] \'Elie CARTAN}, \textit{Les groupes de transformations continus, infinis, simples},
Ann. Sci. \'Ecole Norm. Sup. 26 (1909), 93-161.\\
\textbf{[C-M-S] Quo-Shin CHI, Sergei A. MERKULOV and Lorenz J. SCHWACHHÖFER},
\textit{On the existence of an infinite series of exotic holonomies},
Invent. Math. 126 (1996), 391-411.\\
\textbf{[K-N] Shoshichi KOBAYASHI and Tadashi NAGANO}, \textit{On filtered Lie algebras and 
geometric structures I, II, III}, J. Math. Mech. 13 (1964), 875-908, 14 (1965), 513-521
and 679-706.\\
\textbf{[M-S] Sergei A. MERKULOV and Lorenz J. SCHWACHHÖFER},
\textit{Classification of irreducible holonomies of torsion free connections},
Annals of Mathematics 150 (1999), 77-149 and 1177-1179.\\
\textbf{[W] Joseph A. WOLF},\textit{Complex homogeneous contact manifolds and quaternionic
symmetric spaces}, J. of Math. and Mech. 14 (1965), 1033-1047.

\bigskip
Lionel B\'ERARD BERGERY

Institut \'Elie Cartan à Nancy

Université de Lorraine (FRANCE)

lionel.berard-bergery@univ-lorraine.fr

November $9^{th},\, 2012$
\end{document}